% Format Light AzTeX (F. Arnault)

\catcode`@=11

               %%%%%%%%%%%%%%%%%%%% Fonts %%%%%%%%%%%%%%%%%%%%%

\font\seventeenrm=cmr17
\font\seventeeni=cmmi12 scaled\magstep2
\font\seventeensy=cmsy10  scaled\magstep3
\font\seventeenex=cmex10  scaled\magstep3
\font\seventeengreek=cmr12 scaled\magstep2
\font\seventeenbf=cmbx12 scaled\magstep2
\font\seventeenib=cmmib10 scaled\magstep3
\font\seventeensyb=cmbsy10 scaled\magstep3
\font\seventeenexb=cmexb10 scaled\magstep3

\font\seventeengreekb=cmbx12 scaled\magstep2

\font\fourteenrm=cmr12   scaled\magstep1
\font\fourteeni=cmmi12   scaled\magstep1
\font\fourteensy=cmsy10  scaled\magstep2
\font\fourteenex=cmex10  scaled\magstep2
\font\fourteengreek=cmr12 scaled\magstep1
\font\fourteenbf=cmbx12  scaled\magstep1
\font\fourteensc=cmcsc10 scaled\magstep2
\font\fourteenib=cmmib10 scaled\magstep2

\font\twelverm=cmr12
\font\twelvei=cmmi12
\font\twelvesy=cmsy10  scaled\magstep1
\font\twelveex=cmex10  scaled\magstep1
\font\twelvebf=cmbx12
\font\twelveit=cmti12
\font\twelvebfit=cmbxti12
\font\twelvesl=cmsl12
\font\twelvebfsl=cmbxsl10 scaled\magstep1
\font\twelvett=cmtt12
\font\twelvebftt=cmbtt10 scaled\magstep1
\font\twelvesc=cmcsc10 scaled\magstep1
\font\twelvebfsc=cmbcsc10 scaled\magstep1
\font\twelveib=cmmib10 scaled\magstep1
\font\twelvesyb=cmbsy10 scaled\magstep1
\font\twelveexb=cmexb10 scaled\magstep1
\font\twelvemsam=msam10 scaled\magstep1
\font\twelvebb=msbm10  scaled\magstep1
\font\twelvegoth=eufm10    scaled\magstep1
\font\twelvescript=eusm10  scaled\magstep1
\font\twelvegreek=cmr12
\font\twelvegreekb=cmbx12

\font\tensc=cmcsc10
\font\tenib=cmmib10
\font\tensyb=cmbsy10
\font\tenexb=cmexb10
\font\tenbfsl=cmbxsl10
\font\tenbfit=cmbxti10
\font\tenbfsc=cmbcsc10
\font\tenbftt=cmbtt10
\font\tenmsam=msam10
\font\tenbb=msbm10
\font\tengoth=eufm10
\font\tenscript=eusm10
\font\tengreek=cmr10
\font\tengreekb=cmbx10

\font\ninerm=cmr9
\font\ninei=cmmi9
\font\ninesy=cmsy9

\font\eightrm=cmr8
\font\eighti=cmmi8
\font\eightsy=cmsy8
\font\eightbf=cmbx8
\font\eightsl=cmsl8
\font\eightit=cmti8
\font\eighttt=cmtt8
\font\eightsc=cmcsc8

\font\sevenmsam=msam7
\font\sevenbb=msbm7
\font\sevengoth=eufm7
\font\sevenscript=eusm7
\font\sevengreek=cmr7

\font\fivemsam=msam5

                   % Standard families of fonts (0,1,2,3) %

\def\seventeenm@th{%
 \textfont0=\seventeenrm \scriptfont0=\fourteenrm \scriptscriptfont0=\twelverm
 \textfont1=\seventeeni  \scriptfont1=\fourteeni  \scriptscriptfont1=\twelvei
 \textfont2=\seventeensy \scriptfont2=\fourteensy \scriptscriptfont2=\twelvesy
 \textfont3=\seventeenex \scriptfont3=\fourteenex \scriptscriptfont3=\twelveex}

\def\seventeenm@thbf{%
 \textfont0=\seventeenbf  \scriptfont0=\fourteenbf \scriptscriptfont0=\twelverm
 \textfont1=\seventeenib  \scriptfont1=\fourteenib  \scriptscriptfont1=\twelvei
 \textfont2=\seventeensyb \scriptfont2=\fourteensy \scriptscriptfont2=\twelvesy
 \textfont3=\seventeenexb \scriptfont3=\fourteenex \scriptscriptfont3=\twelveex}

\def\fourteenm@th{%
 \textfont0=\fourteenrm  \scriptfont0=\tenrm  \scriptscriptfont0=\sevenrm
 \textfont1=\fourteeni   \scriptfont1=\teni   \scriptscriptfont1=\seveni
 \textfont2=\fourteensy  \scriptfont2=\tensy  \scriptscriptfont2=\sevensy
 \textfont3=\fourteenex  \scriptfont3=\tenex  \scriptscriptfont3=\tenex}

\def\twelvem@th{%
 \textfont0=\twelverm  \scriptfont0=\ninerm  \scriptscriptfont0=\sevenrm
 \textfont1=\twelvei   \scriptfont1=\ninei   \scriptscriptfont1=\seveni
 \textfont2=\twelvesy  \scriptfont2=\ninesy  \scriptscriptfont2=\sevensy
 \textfont3=\twelveex  \scriptfont3=\tenex   \scriptscriptfont3=\tenex}

\def\twelvem@thbf{%
 \textfont0=\twelvebf  \scriptfont0=\tenbf  \scriptscriptfont0=\sevenbf
 \textfont1=\twelveib  \scriptfont1=\tenib  \scriptscriptfont1=\seveni
 \textfont2=\twelvesyb \scriptfont2=\ninesy  \scriptscriptfont2=\sevensy
 \textfont3=\twelveexb \scriptfont3=\tenex   \scriptscriptfont3=\tenex}

\def\tenm@th{%
 \textfont0=\tenrm  \scriptfont0=\sevenrm  \scriptscriptfont0=\fiverm
 \textfont1=\teni   \scriptfont1=\seveni   \scriptscriptfont1=\fivei
 \textfont2=\tensy  \scriptfont2=\sevensy  \scriptscriptfont2=\fivesy
 \textfont3=\tenex  \scriptfont3=\tenex    \scriptscriptfont3=\tenex}

\def\tenm@thbf{%
 \textfont0=\tenbf  \scriptfont0=\sevenrm  \scriptscriptfont0=\fiverm
 \textfont1=\tenib  \scriptfont1=\seveni   \scriptscriptfont1=\fivei
 \textfont2=\tensyb \scriptfont2=\sevensy  \scriptscriptfont2=\fivesy
 \textfont3=\tenexb \scriptfont3=\tenex    \scriptscriptfont3=\tenex}

\def\eightm@th{%
 \textfont0=\eightrm  \scriptfont0=\fiverm  \scriptscriptfont0=\fiverm
 \textfont1=\eighti   \scriptfont1=\fivei   \scriptscriptfont1=\fivei
 \textfont2=\eightsy  \scriptfont2=\fivesy  \scriptscriptfont2=\fivesy
 \textfont3=\tenex    \scriptfont3=\tenex   \scriptscriptfont3=\tenex}

\def\sevenm@th{%
 \textfont0=\sevenrm  \scriptfont0=\fiverm  \scriptscriptfont0=\fiverm
 \textfont1=\seveni   \scriptfont1=\fivei   \scriptscriptfont1=\fivei
 \textfont2=\sevensy  \scriptfont2=\fivesy  \scriptscriptfont2=\fivesy
 \textfont3=\tenex    \scriptfont3=\tenex   \scriptscriptfont3=\tenex}

%\def\fivem@th{%
% \textfont0=\fiverm   \scriptfont0=\fiverm  \scriptscriptfont0=\fiverm
% \textfont1=\fivei    \scriptfont1=\fivei   \scriptscriptfont1=\fivei
% \textfont2=\fivesy   \scriptfont2=\fivesy  \scriptscriptfont2=\fivesy
% \textfont3=\tenex    \scriptfont3=\tenex   \scriptscriptfont3=\tenex}

                             % font families %

% These families are known by Plain :
% \fam0, \fam1 (\mit, \oldstyle), \fam2 (\cal, symbols), \fam3 (extended sy.)
% \fam\itfam = \fam4 (only \textfont)
% \fam\slfam = \fam5 (only \textfont)
% \fam\bffam = \fam6
% \fam\ttfam = \fam7 (only \textfont)

\newfam\greekf@m
\newfam\msamf@m
\newfam\bbf@m
\newfam\gothf@m
\newfam\scriptf@m

		     % Greek uppercase letters %

\def\Gamma{{\fam\greekf@m\mathchar"7800}}
\def\Delta{{\fam\greekf@m\mathchar"7801}}
\def\Theta{{\fam\greekf@m\mathchar"7802}}
\def\Lambda{{\fam\greekf@m\mathchar"7803}}
\def\Xi{{\fam\greekf@m\mathchar"7804}}
\def\Pi{{\fam\greekf@m\mathchar"7805}}
\def\Sigma{{\fam\greekf@m\mathchar"7806}}
\def\Upsilon{{\fam\greekf@m\mathchar"7807}}
\def\Phi{{\fam\greekf@m\mathchar"7808}}
\def\Psi{{\fam\greekf@m\mathchar"7809}}
\def\Omega{{\fam\greekf@m\mathchar"780A}}

                              % Character sizes %

\def\seventeenpoint{%
  \def\lf{%
     \textfont\greekf@m=\seventeengreek
     \scriptfont\greekf@m=\fourteengreek
     \def\rm{\seventeenm@th\fam0\seventeenrm}%
     \def\it{\sevenit}%
     \def\sl{\sevensl}%
     \def\tt{\seventt}%
     \def\sc{\sevensc}%
     \rm}%
  \def\bf{%
     \textfont\greekf@m=\seventeengreekb
     \scriptfont\greekf@m=\fourteengreek
     %\scriptscriptfont\greekf@m=\twelvegreekb
     \def\rm{\seventeenm@thbf\fam0\seventeenbf}%
     \rm}%
  \normalbaselineskip=20pt\normalbaselines
  \lf}

\def\fourteenpoint{%
  \def\lf{\def\rm{\fourteenm@th\fam0\fourteenrm}\rm}%
  \def\bf{\def\rm{\fourteenbf}\fam0\rm}%
  \def\sc{\fourteensc}%
  \normalbaselineskip=17pt\normalbaselines
  \lf}

\def\twelvepoint{%
  \textfont\msamf@m=\twelvemsam
  \scriptfont\msamf@m=\tenmsam
  %\scriptscriptfont\msamf@m=\sevenmsam
  \textfont\bbf@m=\twelvebb
  \scriptfont\bbf@m=\tenbb
  \scriptscriptfont\bbf@m=\sevenbb
  \textfont\gothf@m=\twelvegoth
  \scriptfont\gothf@m=\tengoth
  \scriptscriptfont\gothf@m=\sevengoth
  \textfont\scriptf@m=\twelvescript
  \scriptfont\scriptf@m=\tenscript
  \scriptscriptfont\scriptf@m=\sevenscript
  \def\lf{%
      \textfont\greekf@m=\twelvegreek
      \scriptfont\greekf@m=\tengreek
      \scriptscriptfont\greekf@m=\sevengreek
      \def\rm{\twelvem@th\fam0\twelverm}%
      \def\it{\twelveit}%
      \def\sl{\twelvesl}%
      \def\tt{\twelvett}%
      \def\sc{\twelvesc}%
      \rm}%
  \def\bf{%
      \textfont\greekf@m=\twelvegreekb
      \scriptfont\greekf@m=\tengreekb
      %\scriptscriptfont\greekf@m=\sevengreekb
      \def\rm{\twelvem@thbf\fam0\twelvebf}%
      \def\it{\twelvebfit}%
      \def\sl{\twelvebfsl}%
      \def\tt{\twelvebftt}%
      \def\sc{\twelvebfsc}%
      \rm}%
  \def\msam{\fam\msamf@m\twelvemsam}%
  \def\bb{\fam\bbf@m\twelvebb}%
  \def\goth{\fam\gothf@m\twelvegoth}%
  \def\script{\fam\scriptf@m\twelvescript}%
  \normalbaselineskip=14pt\normalbaselines
  \lf}

\def\tenpoint{%
  \textfont\msamf@m=\tenmsam
  \scriptfont\msamf@m=\sevenmsam
  \scriptscriptfont\msamf@m=\fivemsam
  \textfont\bbf@m=\tenbb
  \scriptfont\bbf@m=\sevenbb
  %\scriptscriptfont\bbf@m=\fivebb
  \textfont\gothf@m=\tengoth
  \scriptfont\gothf@m=\sevengoth
  %\scriptscriptfont\gothf@m=\fivegoth
  \textfont\scriptf@m=\tenscript
  \scriptfont\scriptf@m=\sevenscript
  %\scriptscriptfont\scriptf@m=\fivescript
  \def\lf{%
      \textfont\greekf@m=\tengreek
      \scriptfont\greekf@m=\sevengreek
      %\scriptscriptfont\greekf@m=\fivegreek
      \def\rm{\tenm@th\fam0\tenrm}%
      \def\it{\tenit}%
      \def\sl{\tensl}%
      \def\tt{\fam\ttfam\tentt}%
      \def\sc{\tensc}%
      \rm}%
  \def\bf{%
      \textfont\greekf@m=\tengreekb
      %\scriptfont\greekf@m=\sevengreekb
      %\scriptscriptfont\greekf@m=\fivegreekb
      \def\rm{\tenm@thbf\fam0\tenbf}%
      \def\it{\tenbfit}%
      \def\sl{\tenbfsl}%
      \def\tt{\fam\ttfam\tenbftt}%
      \def\sc{\tenbfsc}%
      % will redefine \goth,\script ?
      \rm}%
  \def\msam{\fam\msamf@m\tenmsam}%
  \def\bb{\fam\bbf@m\tenbb}%
  \def\goth{\fam\gothf@m\tengoth}%
  \def\script{\fam\scriptf@m\tenscript}%
  \normalbaselineskip=12pt\normalbaselines
  \lf}

\def\eightpoint{%
  \def\lf{\def\rm{\eightm@th\fam0\eightrm}\rm}%
  \def\bf{\def\rm{\eightbf}\rm}%
  \def\it{\eightit}%
  \def\sl{\eightsl}%
  \def\tt{\eighttt}%
  \def\sc{\eightsc}%
  \normalbaselineskip=10pt\normalbaselines
  \setbox\strutbox=\hbox{\vrule height7pt depth3pt width0pt}%
  \lf}

\def\sevenpoint{%
  \def\lf{\def\rm{\sevenm@th\fam0\sevenrm}\rm}%
  \def\bf{\def\rm{\sevenbf}\rm}%
  \def\it{\sevenit}%
  \def\sl{\sevensl}%
  \def\tt{\eighttt}% ersatz
  \normalbaselineskip=8pt\normalbaselines
  \lf}

\tenpoint

\def\mathbf#1{\hbox{\bf$#1$}}

               %%%%%%%%%%%%%%% External files %%%%%%%%%%%%%%%%%

\newif\ifreffile          \reffilefalse
\newread\reffile
\newread\pagenosfile
\newwrite\reffile

\outer\def\openreffile{%
   \re@dreffile
   \immediate\write16{Writing references to \jobname.ref}%
   \immediate\openout\reffile=\jobname.ref
   \reffiletrue}

\def\re@dreffile{%
     \openin\reffile=\jobname.ref
     \ifeof\reffile
        \closein\reffile
        \immediate\write16{*** Not found \jobname.ref ***}%
     \else\closein\reffile
        \immediate\write16{Reading references from \jobname.ref}%
     \let\prefixforref=\relax
     % for some Plain internal cs (dimens,...) ar'nt expanded :
     \catcode`@=11
     \input\jobname.ref
     \catcode`@=12
     \fi}

\newread\xreffile
\outer\def\readref#1{%
   \openin\xreffile=#1.ref
   \ifeof\xreffile
       \closein\xreffile
       \immediate\write16{*** Not found #1.ref ***}%
   \else\closein\xreffile
        \immediate\write16{Reading references from #1.ref}%
        \def\prefixforref{\xchapno.}%
        \catcode`@=11  % for some Plain internal cs (dimens,...) ar'nt expanded
        \input#1.ref
        \catcode`@=12
        \let\prefixforref=\relax
   \fi
   }

               %%%%%%%%%%%%%%%%%% Page setting %%%%%%%%%%%%%%%%

%\hoffset=-3mm    % beurk !
%\voffset=8mm    % (felicie)

%\hsize=6.5in
%\vsize=8.9in
\def\landscape{{\let\tmpsize=\hsize
                \global\let\hsize=\vsize
		\global\let\vsize=\tmpsize}}

\def\t@define{$\spadesuit$}

\def\getp@geno{\global\pageno=\pagenos\chapterno}

\newif\ifshowref    \showreffalse
\let\innerl@bel=\relax
\def\label#1{\def\innerl@bel{#1}\par}  % \par is useful after \subsection
                                       % but I don't now why it works...

                          % headlines and footlines %

\def\frenchmonth{\ifcase\month\or janvier\or f\'evrier\or mars\or avril%
       \or mai\or juin\or juillet\or ao\^ut%
       \or septembre\or octobre\or novembre\or d\'ecembre\fi}
\def\englishmonth{\ifcase\month\or January\or February\or March\or April%
       \or May\or June\or July\or August%
       \or September\or October\or November\or December\fi}
\def\frenchdate{\the\day\ \frenchmonth\ \the\year}
\def\englishdate{\englishmonth\ \the\day, \the\year}
\let\date=\frenchdate
\newif\ifshowdate       \showdatetrue
\def\author{F.~Arnault}
\newif\iffirstpage      \firstpagetrue
\newtoks\firstheadline
\newtoks\rightheadline
\newtoks\leftheadline
\newtoks\firstfootline
\newtoks\otherfootline

\headline={%
   \iffirstpage\the\firstheadline
   \else\ifodd\pageno\the\rightheadline
        \else\the\leftheadline\fi
   \fi}
\firstheadline={\hfil}
\rightheadline={\tenpoint\sl\hfil\firstmark\hfil}
\leftheadline={\tenpoint\sl\hfil\firstmark\hfil}

\footline={%
   \iffirstpage\the\firstfootline\global\firstpagefalse
   \else\the\otherfootline
   \fi}
\firstfootline={\ifshowdate \hfil\sevenrm\author\quad---\quad\date\hfil 
                \else\footfolio \fi}
\otherfootline={\footfolio}
\def\footfolio{\hfil\tenrm --- \folio\ ---\hfil}

			% Chapters, titles, subtitles %

\newskip\titleskipamount
%\titleskipamount=36pt plus4pt minus4pt  % 3\bigskip
\titleskipamount=24pt plus4pt minus4pt  % 2\bigskip
\def\chapterno{\t@define}
\def\ch@pterm@rk{\relax}

\outer\def\chapter#1#2 #3\par{%  %% number, name, title.
   \global\edef\chapterno{#1}%
   \ifnum#1>-1
       \global\edef\ch@pterm@rk{\chapterno.}%
   \fi
   \re@dp@genos
   \getp@geno
   \ifreffile
      \immediate\write\reffile{\def\string\xchapno{#1}}%
      \immediate\write\reffile{\def\string\chap#2{#1}}%
   \fi
   \message{#2 -> chapitre #1 : #3. }%
   % headline of pages 2, 3, ... if no section begins :
   {\def\cr{\relax}\mark{#3}}%
   \leftheadline={\tenpoint\sl\hfil\ch@pterm@rk{\def\cr{\relax} #3}\hfil}
   \ifnum#1>-1
     \leftline{\twelvepoint\bf Chapitre #1}%
     \bigskip
   \fi
   {\twentyonepointb\halign{\centerline{##}\cr#3\crcr}}%
   \ifreffile
      \write\reffile{\string\t@cchapter{\folio}{#1}{#3}}%
   \fi
   \bigskip\bigskip\bigskip
   }

\outer\def\title#1\par{%
   % headline of pages 2, 3, ... if no section begins : 
   {\def\cr{\relax}\mark{#1}}%
   {\seventeenpoint\bf\halign{\centerline{##}\cr#1\crcr}}%
   \ifreffile
      \immediate\write\reffile{\def\string\xchapno{\relax}}%
      \immediate\write\reffile{\string\t@ctitle{#1}}%
   \fi
   \vskip\titleskipamount}

\outer\def\subtitle#1\par{%
   \vskip-\titleskipamount
   {\fourteenpoint\halign{\centerline{##}\cr#1\crcr}}%
   \vskip\titleskipamount}

\parskip=2pt plus 1pt minus 0.5pt   %% Plain : 0pt plus 1pt

				 % Abstracts %

\newdimen\abstractmarginswidth
\abstractmarginswidth=30pt
\long\def\abstract#1.#2\endabstract{%
  \begingroup
  \parindent=\abstractmarginswidth
  \narrower\noindent\eightpoint
  {\bf#1.\enspace}%
  #2\par
  \endgroup}
\def\endabstract{%
 \errhelp={I have encountered a end-of-abstract mark, outside of any abstract}%
 \errmessage{\string\endabstract\space must match an \string\abstract}}

			   % Sections, subsections %

\newif\ifromansection \romansectionfalse
\def\stringsectionnumber#1{\ifromansection\uppercase\expandafter{\romannumeral#1}\else\the#1\fi}

\newcount\sectionnumber
\outer\def\section#1 \par{% % The first <cr> is converted to <sp>_10, knuth p47
    \global\advance\sectionnumber by1
    \global\statnumber=0
    %\bigbreak but always with penalty
    \ifdim\lastskip<\bigskipamount \removelastskip\bigskip \fi
    \penalty-400
    \ifx\innerl@bel\relax\relax
    \else
      \expandafter\xdef\csname sec\innerl@bel\endcsname{\the\sectionnumber}%
      \message{\innerl@bel -> section \the\sectionnumber, }%
      \ifreffile
         \immediate\write\reffile{\edef\string\sec\innerl@bel{%
                          \string\prefixforref\stringsectionnumber{\sectionnumber}}}%
      \fi
    \fi
    \leftline{%
          \ifshowref
            \ifx\innerl@bel\relax\relax  \else\llap{\fiverm\innerl@bel\ }\fi
          \fi
          \twelvepoint\bf\stringsectionnumber{\sectionnumber}. #1}%
    \mark{\ch@pterm@rk\stringsectionnumber{\sectionnumber}. #1}%
    \ifreffile
       \write\reffile{\string\t@csecti@n{\folio}%
                                        \string{#1\string}}%
    \fi
    \nobreak\medskip\noindent
    \let\innerl@bel=\relax
    }

\outer\def\subsection#1 \par{%
    \ifhmode  %% Probably, the \subsection just follows a \section (\noindent)
       \immediate\write-1{%
             I guess that the subsection just follows
             section \the\sectionnumber's title.}%
    \else
        %\bigbreak but always with penalty
        \ifdim\lastskip<\bigskipamount \removelastskip\bigskip \fi
        \penalty-300 % changed from -200
    \fi
    \leftline{\twelvepoint#1}%
    \nobreak\medskip\noindent
    }

\outer\def\subsubsection#1 \par{%
    \ifhmode  %% Probably, the \subsection just follows a \section (\noindent)
       \immediate\write-1{%
             I guess that the subsubsection just follows a section title.}%
    \else
        %\bigbreak but always with penalty
        \ifdim\lastskip<\bigskipamount \removelastskip\bigskip \fi
        \penalty-250 % changed from -200
    \fi
    \leftline{#1}%
    \nobreak\smallskip\noindent
    }

                        % Statements, equations, ... %

\newcount\statnumber

\def\innerstat#1#2{%
    \global\advance\statnumber by1
    \ifx\innerl@bel\relax\relax
    \else
       \expandafter\xdef\csname ref\innerl@bel\endcsname
          {\stringsectionnumber{\sectionnumber}.\the\statnumber}%
       \message{\innerl@bel -> \csname ref\innerl@bel\endcsname, }%
       \ifreffile\immediate\write\reffile{%
             \edef\string\ref\innerl@bel{%
             \string\prefixforref\stringsectionnumber{\sectionnumber}.\the\statnumber}}%
       \fi
    \fi
    \medbreak\noindent
    \ifshowref
       \ifx\innerl@bel\relax\relax \else\llap{\fiverm\innerl@bel\ }\fi
    \fi
    \alphaenumreset
    {\bf\stringsectionnumber\sectionnumber.\the\statnumber. --- #1.\enspace}%
    %{\bf\romansection\the\sectionnumber.\the\statnumber. --- #1.\enspace}%
    {#2\par}%
    \let\innerl@bel=\relax
    \ifdim\lastskip<\medskipamount \removelastskip\penalty55\medskip\fi
    }

\outer\def\stat #1. #2\par{\innerstat{#1}{\sl#2}}
\outer\def\statrm #1. #2\par{\innerstat{#1}{#2}} 
\outer\def\remark #1. #2\par{%
    \smallbreak\noindent
    {\bf#1.\enspace}%
    {#2\par}%
    \ifdim\lastskip<\medskipamount \removelastskip\penalty55\medskip\fi
    }

\newcount\eqnumber
\def\eqdef#1{%
    \global\advance\eqnumber by1
    \expandafter\xdef\csname eq#1\endcsname{\the\eqnumber}%
    \message{#1 -> (\csname eq#1\endcsname), }%
    \ifreffile\immediate\write\reffile{%
          \edef\string\eq#1{\string\prefixforref\the\eqnumber}}%
    \fi
    \eqno\hbox{%
       \lf(\the\eqnumber)%
       \ifshowref\rlap{\fiverm\ #1}\fi}}

% for use in \eqalignno :
\def\eqaligndef#1{% 
    \global\advance\eqnumber by1
    \expandafter\xdef\csname eq#1\endcsname{\the\eqnumber}%
    \message{#1 -> (\csname eq#1\endcsname), }%
    \ifreffile\immediate\write\reffile{\def\string\eq#1{\the\eqnumber}}%
    \fi
    \lf(\the\eqnumber)%
    \ifshowref\rlap{\fiverm\ #1}\fi}

% enumerations
\newcount\enumnumber

\def\enum{%
   \advance\enumnumber by1
   \par\noindent(\the\enumnumber) }

\def\alphaenumreset{\enumnumber="60 }
\def\alphaenum{%
  \advance\enumnumber by1
  \par\noindent(\char\enumnumber) }
\def\rawalphaenum{%
  \advance\enumnumber by1
  (\char\enumnumber) }

			      % Figures %

\newbox\figb@x
\newcount\fignumber
\newdimen\figb@xwidth

\def\figure#1. #2#3{%
  \global\advance\fignumber by1
  \ifx\innerl@bel\relax\relax
  \else
     \expandafter\xdef\csname fig\innerl@bel\endcsname{\the\fignumber}%
     \message{\innerl@bel -> #1 \csname fig\innerl@bel\endcsname, }%
     \ifreffile\immediate\write\reffile{%
         \edef\string\fig\innerl@bel{%
         \string\prefixforref\the\fignumber}}%
     \fi
  \fi
  \setbox\figb@x=\vbox{#3}%
  \figb@xwidth=\wd\figb@x
  \setbox\figb@x=\vbox{%
       \box\figb@x
       \smallskip
       \hbox to\figb@xwidth{%
          \ifshowref
             \ifx\innerl@bel\relax\relax \else\llap{\fiverm\innerl@bel\ }\fi
          \fi
          \hss{\bf #1\ \the\fignumber.} #2\hss}}%
  $$%
    \box\figb@x
  $$%
  \let\innerl@bel=\relax}

                              % Bibliographies %

\newbox\bibb@x
\newcount\bibchosen
\newcount\bibnumber
\outer\def\showbib{\unvbox\bibb@x}

\def\bibdef#1#2#3#4{%
    \global\advance\bibnumber by1
    \expandafter\xdef\csname bib#1\endcsname{\the\bibnumber}%
    \global\setbox\bibb@x=\vbox{%
        \unvbox\bibb@x
        \hang\textindent{[\the\bibnumber]}%
        \strut{\sc#2 :\ \it#3.\/\ \ \rm#4.}%
        \par\medskip}}

\def\bibunlo@ded{\t@define}

\def\bibchoose#1{%
    \global\advance\bibchosen by1
    \expandafter\def\csname bib#1\endcsname{\t@define}}

\let\ifbibchooseall=\iffalse
\def\bibchooseall{\let\ifbibchooseall=\iftrue}

\outer\def\checkbib{\ifbibchooseall\relax
         \else\ifnum\bibchosen=\bibnumber\relax
              \else\immediate\write16{Warning : Some bibliographic 
                                      references remain undefined.}%
              \fi
         \fi}

\outer\def\bibitem#1#2#3#4{%
   \ifbibchooseall\bibdef{#1}{#2}{#3}{#4}%
   \else
       \expandafter\ifx\csname bib#1\endcsname\relax \relax %
       \else
           \expandafter\ifx\csname bib#1\endcsname\bibunlo@ded
               \bibdef{#1}{#2}{#3}{#4}%
           \else\immediate\write16{Warning : attempt to redefine
                                   bibliographic reference #1.}%
           \fi
       \fi
   \fi}

                      % Reference substitution %

% These must be defined AFTER \chapter, \eqdef, \figure...

\def\chapref#1{%
    \expandafter\ifx\csname chap#1\endcsname\relax \t@define
          \immediate\write16{ *** Undefined chapter reference : #1 ***}%
    \else\csname chap#1\endcsname
    \fi}

\def\secref#1{%
    \expandafter\ifx\csname sec#1\endcsname\relax \t@define
          \immediate\write16{ *** Undefined section reference : #1 ***}%
    \else\csname sec#1\endcsname
    \fi}

\def\ref#1{%
    \expandafter\ifx\csname ref#1\endcsname\relax \t@define
          \immediate\write16{ *** Undefined statement reference : #1 ***}%
    \else\csname ref#1\endcsname
    \fi}

\def\eqref#1{%
    \expandafter\ifx\csname eq#1\endcsname\relax \t@define
          \immediate\write16{ *** Undefined equation reference : #1 ***}%
    \else(\csname eq#1\endcsname)%
    \fi}

\def\figref#1{%
    \expandafter\ifx\csname fig#1\endcsname\relax \t@define
          \immediate\write16{ *** Undefined figure reference : #1 ***}%
    \else\csname fig#1\endcsname
    \fi}

\def\bib#1{%
   \expandafter\ifx\csname bib#1\endcsname\relax
      [\t@define]%
      \immediate\write16{ *** Unknown bibliographic reference : #1 ***}%
   \else
      \expandafter\ifx\csname bib#1\endcsname\bibunlo@ded
          \immediate\write16{ *** Unloaded bibliographic reference : #1 ***}%
      \fi
      [\csname bib#1\endcsname]%
   \fi}

               %%%%%%%%%%%%%%% Tables of Contents %%%%%%%%%%%%%

\newif\iftoc   \tocfalse
\newbox\t@cb@x
\newcount\t@csecn@
\newif\ifnarrowtoc  \narrowtoctrue

\def\t@cchapter#1#2#3{%  page, chapno, titre.
  \iftoc{%
     \global\narrowtocfalse
     \def\cr{\relax}%
     \advance\hsize by -\pagecolwidth
     \global\t@csecn@=0
     \parindent=0pt
     \global\setbox\t@cb@x=\vbox{%
         \unvbox\t@cb@x
         \bigskip\bigskip\goodbreak
         \ifnum#2>-1
           \centerline{\bf Chapitre #2}%
           \smallskip\nobreak
         \fi
         \strut\blacksquare\ \bf#3\page{#1}\kern-\pagecolwidth}%
         \smallskip\nobreak}%
  \else\relax
  \fi}

\def\t@ctitle#1{%
  \iftoc{%
     \global\narrowtoctrue
     \parindent=\abstractmarginswidth
     \narrower
     \def\cr{\relax}%
     \advance\hsize by -\pagecolwidth
     \global\t@csecn@=0
     \parindent=0pt
     \global\setbox\t@cb@x=\vbox{%
         \unvbox\t@cb@x
         \medskip
         \strut\blacksquare\ \bf#1}%
         \medskip\nobreak}%
  \else\relax
  \fi}

\def\t@csecti@n#1#2{%
  \iftoc{%
     \advance\hsize by -\pagecolwidth
     \global\advance\t@csecn@ by1
     \parindent=0pt
     \global\setbox\t@cb@x=\vbox{%
         \unvbox\t@cb@x
             \strut
             \ifnarrowtoc\kern\abstractmarginswidth\fi
             \the\t@csecn@. --- #2\page{#1}%
             \ifnarrowtoc\kern\abstractmarginswidth\fi
             \kern-\pagecolwidth}%
          \smallskip}%
  \else\relax
  \fi}

\def\showtoc{%
  \iftoc\unvbox\t@cb@x
  \else
    \immediate\write16{*** TOC is void because you hav'nt said
                       \string\toctrue. ***}%
  \fi}

               %%%%%%%%%%%%%%%%%% \everyjob %%%%%%%%%%%%%%%%%%%

\def\pagenos#1{1\write16{Using default pageno : 1}}

\def\re@dp@genos{%
     \openin\pagenosfile=pagenos.tex
     \ifeof\pagenosfile\closein\pagenosfile
     \else\closein\pagenosfile
        \immediate\write16{Reading pagenos.tex}%
     % pagenos.tex is assumed to contain the following definition
     % \def\pagenos#1{\ifcase#1 n1 \or n2 \or ... \else 1 \fi}%
     \input pagenos.tex
     \fi}

\everyjob={%
     \immediate\write16{Format Light AzTeX (F. Arnault), version \fmtversion.}
     \setbox\t@cb@x=\vbox{}%
     \setbox\bibb@x=\vbox{}%
     \tracingstats=1}

               %%%%%%%%%%%% Verbatim and Listings %%%%%%%%%%%%%

% see TeXbook, pages 380 -->

% The following cannot be included in a \def because
% the assignment can be parsed only if \obeyspace has expanded.
{\obeyspaces\global\let =\ }  % Active space will have fixed width
% In Plain, active space is \space.  But ' ' is usually not active.
{\def\tmpminus{-\relax}\catcode`-=\active \xdef-{\tmpminus}}
{\def\tmpgreater{>\relax}\catcode`>=\active \xdef>{\tmpgreater}}
{\def\tmplessthan{<\relax}\catcode`<=\active \xdef<{\tmplessthan}}

% Assigns all (special) characters to letter category
\def\unc@tcodespeci@ls{\def\do##1{\catcode`##1=12}\dospecials}

% Verbatim settings will be on before parameter is read
\outer\def\verbatim{\par\begingroup\setupverb@tim\doverbatim}

\def\setupverb@tim{%
  \parindent=0mm
  \def\par{\leavevmode\endgraf}%  % Usual \par is \endgraf
  % Spaces and - will be active :
  \obeylines\unc@tcodespeci@ls\catcode`-=\active\catcode`>=\active
       \catcode`<=\active\obeyspaces\tt}

% The only way to escape is \endverbatim at the end of a line.
% \endverbatim will be seen as a sequence of letters, not a \cs
{\catcode`\|=0 \catcode`\\=12  % Switch escape char to |
 |obeylines|gdef|doverbatim^^M#1\endverbatim
 {|medskip|hrule|nobreak#1|smallskip|nobreak|hrule|medskip|endgroup|noindent}}

\def\endverbatim{%
  \errhelp={I have encountered an end-verbatim mark whereas I was not in
           verbatim mode}%
  \errmessage{\string\endverbatim\space must match a \string\verbati m}}

{\catcode`\|=\active\obeylines%
\gdef\|{\begingroup\catcode`\|=\active\setupverb@tim\let^^M=\ \let|=\endgroup}}

\newcount\l@st@ngl@ne
\outer\def\listing#1#2\par{%
   \par\begingroup\setupverb@tim
   \everypar{\advance\l@st@ngl@ne by1\llap{\sevenrm\the\l@st@ngl@ne\quad}}%
   \medskip
   \line{\leaders\hrule\hfil\quad\sl#1\quad\leaders\hrule\hfil}\nobreak
   \input#2
   \nobreak\smallskip\nobreak\hrule\medskip
   \endgroup}

% Cannot make listing breaking work.
%\newinsert\verbinsert
%\skip\verbinsert=\z@skip
%\count\verbinsert=1000
%\dimen\verbinsert=10pt % only a limit
%\newinsert\lastverbinsert
%\skip\lastverbinsert=\z@skip
%\count\lastverbinsert=1000
%\dimen\lastverbinsert=10pt % only a limit
%\output{\verboutput}
%\def\verboutput{\shipout\vbox{\makeheadline\copy\lastverbinsert\pagebody\copy\verbinsert\makefootline}%
%  \advancepageno
%  \ifnum\outputpenalty>-\@MM \else\dosupereject\fi}

                                 % Hidden parts %

\newif\ifshowhidden  \showhiddenfalse

\def\endhidden{\ifshowhidden\relax\else
  \errhelp={I have encountered a end-of-hidden-part mark, 
            without matching hidden-part mark.}%
  \errmessage{\string\endhidden\space must match a \string\hidde n}\fi}

\outer\def\hidden{\ifshowhidden\relax
                  \else\par\begingroup\s@tuphidd@n\dohidden\fi}

\def\s@tuphidd@n{\obeylines\unc@tcodespeci@ls}  % I don't understand.

{\catcode`\|=0 \catcode`\\=12  % Switch escape char to |
 |obeylines|gdef|dohidden#1\endhidden{|endgroup|fi}}  % why |fi ???
% \endhidden will be seen as a sequence of letters, not a \cs

               %%%%%%%%%%%%%%%%% Initials %%%%%%%%%%%%%%%%%%%%

\font\initialfont=cmbx12 scaled\magstep4
\newbox\initialb@x
\def\initial#1{%
   \setbox\initialb@x=\hbox{\initialfont#1\hskip2pt}%
   \hang\hangafter=-2\hangindent=\wd\initialb@x
   \setbox\initialb@x=\hbox{\hskip-\wd\initialb@x\box\initialb@x}%
   \noindent
   \smash{\lower12pt\box\initialb@x}%
   }
 
\catcode`@=12

               %%%%%%%%%%%%%%%%%%% Divers %%%%%%%%%%%%%%%%%%%%%

\newdimen\pagecolwidth  \pagecolwidth=8mm
\def\page#1{%
       \quad
       %\leaders\hbox to 3mm{\hfil.\hfil}\hfill
       \leaders\hrule height0.1pt\hfill
       %\quad\hbox{#1}}
       \hbox to\pagecolwidth{\hfill#1}}
       
\def\qed{\kern 4pt\penalty500\null\hfill\square\par}

\def\square{%
   \hbox{%
      \vrule
      \vbox to 6pt{\hrule width 4pt\vfill\hrule}%
      \vrule}}

\def\blacksquare{%
  \vbox{%
    \hbox{%
      \kern1pt
      \vrule height5pt depth 0pt width 5pt
      \kern1pt}%
    \kern0.5pt}}

\def\dem{\vskip-4pt\noindent{\sc D\'emonstration --- }\alphaenumreset}
\def\pdem#1-{\vskip-4pt\noindent{\sc#1 }--- \alphaenumreset}

\def\N{{\bb N}}
\def\Q{{\bb Q}}

\def\Z{{\bb Z}}

% See TeXBook p.154.  Category must be 7 to use non-standard \fam (\bb).
\def\notmid{\mathrel{\bb\mathchar"712D}}
  % semi-direct product
  % k with double bar
\def\hbar{{\mathord{\bb\mathchar"717D}}}  % reduced Planck constant (ugly in Plain)
    % 2nd Hebrew letter
   % 3rd Hebrew letter
  % 4th Hebrew letter
                        % \mathchardef\parallel="326B   is in Plain
\def\Bigmid{\mathrel{\Big|}}
% \mathchardef\leq="3214 \let\le=\leq     % Plain definitions 
% \mathchardef\geq="3215 \let\ge=\geq
\def\leq{\mathrel{\msam\mathchar"7136}}
\def\geq{\mathrel{\msam\mathchar"713E}}
  % double head arrow

%\def\lshift{\mathop{\mathchar"021C}}   % \ll in plain
%\def\rshift{\mathop{\mathchar"021D}}   % \gg in plain
% see also files msam.map, msbm.map, amsfndoc.dvi

  % sometimes better than |#1|

\def\Legendre#1#2{\mathchoice{#1\overwithdelims()#2}{(#1/#2)}{}{}}

\def\smallmatrix#1#2#3#4{\bigl({#1\atop#3}{#2\atop#4}\bigr)}

\def\assign{\mathrel{:=}}
\def\spot{\par\noindent$\bullet$ }

% \item{#1} defined in plain is also useful !

\def\\{\par}  % hidden \par

% \def\gcd{\mathop{\rm gcd}}  % yet in Plain
    % NOT in Plain
\def\pgcd{\mathop{\rm pgcd}}

\def\disc{\mathop{\rm disc}\nolimits}

  % \hom is in Plain

% \def\ker{\mathop{\rm ker}\nolimits} % yet in Plain

  % j'hesite

% \def\deg{\mathop{\rm deg}} % yet in Plain
% \def\det{\mathop{\rm det}} % yet in Plain

\def\SL{\mathop{\rm SL}\nolimits}

\def\fmtversion{2.4.0} % 12/7/11

\showdatefalse

\bibdef{Bernard}
 {A.~Bernard}
 {Formes quadratiques binaires et applications cryptographiques}
 {Th\`{e}se de Doctorat, Universit\'{e} de Limoges, 2011}

\bibdef{BuchmannVollmer}
 {J.~Buchmann, U.~Vollmer}
 {Binary Quadratic Forms, an algorithmic approach}
 {Algorithms and Computation in Mathematics, Vol.~20, Springer, 2007}

\bibdef{Buell}
 {D.A.~Buell}
 {Binary Quadratic Forms, Classical theory and modern computations}
 {Springer, 1989}

\bibdef{CastagnosLaguillaumie}
 {G.~Castagnos, F.~Laguillaumie}
 {On the security of cryptosystems with quadratic decryption: the nicest cryptanalysis}
 {Proceedings of Eurocrypt'09}

\bibdef{CLJN}
 {G.~Castagnos, F.~Laguillaumie, A.~Joux, P.Q.~Nguyen}
 {Factoring with quadratic forms: nice cryptanalyses}
 {Proceedings of Asiacrypt'09}

\bibdef{Cox}
 {D.A.~Cox}
 {Primes of the Form $x^2+ny^2$}
 {Wiley-Interscience, John Wiley \& Sons, 1989}

\bibdef{JSW}
 {M.J.~Jacobson Jr, R.~Scheidler, D.~Weimer}
 {An adaptation of the NICE cryptosystem to real quadratic orders}
 {Proceedings of Africacrypt'08, LNCS~5023, 191--208, 2008}

\bibdef{LehmerLehmer}
 {D.H.~Lehmer, E.~Lehmer}
 {A new factorization technique using quadratic forms}
 {Mathematics of Computation~28, 625--635, 1974}

\bibdef{PaulusTakagi} % Nice (imaginary)
 {S.~Paulus, T.~Takagi}
 {A new public-key cryptosystem over a quadratic order with quadratic decryption time}
 {Journal of Cryptology~13(2), 263--272, 2000}

\bibdef{Schoof}
 {R.J.~Schoof}
 {Quadratic Fields and factorization}
 {Computational Methods in Number Theory, part~II, H.W.~Lenstra, R~Tijdeman
  (ed.).  Mathematical Centre Tracts~155, Amsterdam, 235--286, 1982}

\bibdef{SchnorrLenstra}
 {C.P.~Schnorr, H.W.~Lenstra}
 {A Monte Carlo factoring algorithm with linear storage}
 {Mathematics of Computation~167(43), 289--311, 1984}

\bibdef{Shanks}
 {D.~Shanks}
 {Five number-theoretic algorithms}
 {Proceedings of the Second Manitoba Conference on Numerical Mathematics, 51--70, 1972}

\def\SLZ{\SL_2(\Z)}

\title Formes quadratiques \cr
       de discriminants embo\^{\i}t\'{e}s \cr

\vskip-2mm
\centerline{{\twelvepoint Fran\c cois {\sc Arnault}}%
\footnote*{Adresse \'electronique : {\tt arnault@unilim.fr}}}
\medskip
\centerline{\eightpoint Universit\'e de Limoges --- XLIM (UMR CNRS 6172)}
\centerline{\eightpoint 123 avenue Albert Thomas, F-87060 Limoges Cedex, France}
\bigskip\bigskip

\abstract R\'esum\'e.
  Les formes quadratiques binaires ont \'{e}t\'{e} initialement consid\'{e}r\'{e}es par Fermat,
Lagrange, Legendre.  Puis Gauss, dans les {\it Disquisitiones Arithmeticae} publi\'{e}es en 1801,
est le premier \`{a} leur donner un d\'{e}veloppement significatif, avec en particulier la loi de
composition.

  Leurs applications pratiques sont multiples.  Elles fournissent une mani\`{e}re explicite de
manipuler des id\'{e}aux de corps quadratiques.  De nombreux algorithmes de factorisation les
utilisent : \bib{LehmerLehmer}\bib{Schoof}\bib{SchnorrLenstra}\bib{Shanks}.
Elles sont aussi utilis\'{e}es en cryptographie, en particulier pour les syst\`{e}mes {\sc nice}
%\bib{HPT},
\bib{PaulusTakagi} puis \bib{JSW}.  Les syst\`{e}mes de chiffrement {\sc nice} utilisent des formes
quadratiques de discriminants $\pm p$ et $\pm pq^2$ o\`{u} $p$ et~$q$ sont des nombres premiers.
Cet article pr\'{e}cise les liens entre les formes de discriminant $D$ et celles de discriminant
$Df^2$ (avec $f>1$ entier), ce qui est essentiel pour l'analyse de {\sc nice} et de ses attaques
\bib{Bernard}\bib{CastagnosLaguillaumie}\bib{CLJN}.
Il introduit aussi la notion de formes quadratiques semi-\'{e}quivalentes et en explicite plusieurs
caract\'{e}risations, utiles pour l'analyse de ces attaques~\bib{Bernard}.

\endabstract

\section Introduction

  Cette section rappelle quelques d\'{e}finitions et r\'{e}sultats simples sans d\'{e}monstration.  Elle se
limite volontairement \`{a} ce qui est indispensable pour la suite, passant donc sous silence
des pans entiers de la th\'{e}orie, comme la notion de r\'{e}duction et la loi de groupe sur les classes
d'\'{e}quivalence.  Parmi les nombreux ouvrages de r\'{e}f\'{e}rence sur les formes quadratiques, je
mentionne~\bib{BuchmannVollmer}\bib{Buell}\bib{Cox}.

\subsection Formes

  Une {\it forme quadratique binaire} est un polyn\^{o}me homog\`{e}ne \`{a} deux variables
$$
  q(x, y)
  =
  ax^2 + bxy + cy^2.
$$
Le cas qui nous int\'{e}resse, par la richesse de son arithm\'{e}tique, est celui o\`{u} les coefficients
$a$, $b$, $c$ sont entiers. 
Nous utiliserons le terme abr\'{e}g\'{e} de {\it forme} pour forme quadratique binaire \`{a}
coefficients entiers.  On suppose de plus que le {\it discriminant}
$$
  \disc(q)=D\assign b^2-4ac
$$
n'est pas un carr\'{e} parfait (puisque les propri\'{e}t\'{e}s des formes
quadratiques de discriminant~$D$ sont li\'{e}es au corps sous-jacent $\Q(\sqrt D)$).
Si $D>0$, on parle de forme quadratique {\it r\'{e}elle}.  Si~$D<0$, on parle de forme quadratique
{\it imaginaire}, et on ne s'int\'{e}resse dans ce cas qu'aux formes d\'{e}finies positives (i.e. telles
que~$a>0$).
Nous noterons souvent une forme quadratique en listant ses coefficients $q=(a,b,c)$.
Une forme quadratique $(a,b,c)$ est dite {\it primitive} si $\pgcd(a,b,c)=1$.

\subsection Discriminants

\label{funddisc}
\statrm D\'{e}finitions.  J'appelle {\it discriminant} un entier non carr\'{e} parfait et congru \`{a}~0
ou~1 modulo~4.  On appelle {\it discriminant fondamental} un entier non carr\'{e}
parfait v\'{e}rifiant
$$
  \left\{\eqalign{
     &D \equiv 1 \ \hbox{modulo 4} \cr
     &\hbox{$D$ sans facteur carr\'{e}} \cr}
  \right.
  \qquad\hbox{ou}\qquad
  \left\{\eqalign{
     &D \equiv 0 \ \hbox{modulo 4} \cr
     &\hbox{$D/4\equiv2$ ou~3 modulo 4 et sans facteur carr\'{e}} \cr}
  \right.
$$

  Un discriminant de forme quadratique est donc un discriminant au sens
de la d\'{e}finition~\ref{funddisc}.  Inversement, tout non carr\'{e} parfait
congru \`{a}~0 ou~1 modulo~4 est un discriminant de forme quadratique (par exemple
$x^2-(D/4)y^2$ si $D\equiv0$ modulo~4, et $x^2+xy+((1-D)/4)y^2$ si $D\equiv1$
modulo~4).

  On montre facilement qu'un discriminant $D$ est fondamental si et seulement si il n'existe pas de
discriminant de la forme $D/f^2$ avec $f>1$ entier.  On montre aussi que $D$ est fondamental si et
seulement si toutes les formes quadratiques de discriminant~$D$ sont primitives.

\label{normprimetof}
\stat Lemme.  Soient $q$ une forme quadratique primitive et $f$ un entier non nul.  Alors $q$ est
\'{e}quivalente \`{a} une forme $q'=ax^2+bxy+cy^2$ avec $\pgcd(a,f)=1$.  \qed

\subsection Action, \'{e}quivalence

  Soient $q$ une forme et $M=\smallmatrix prst$ une matrice (\`{a} coefficients r\'{e}els, et non
n\'{e}cessairement de d\'{e}terminant~1).  On note $q\cdot M$ le trin\^{o}me $q'$ donn\'{e} par
$$
  q'(x,y)
  =
  q(px + ry, sx + ty)
  =
  q(p, s) x^2  +  \big(q(p + r, s + t) - q(p, s) - q(r, t)\big) xy  +  q(r, t) y^2.
  \eqdef{coeffsequivform}
$$
On obtient ainsi une action \`{a} droite sur les formes.

\label{discaction}
\stat Proposition.  Si $q'=q\cdot M$ alors $\disc q'=(\det M)^2\disc q$.  \qed
  
  Lorsque la matrice $M$ est \`{a} coefficients entiers de de d\'{e}terminant~1, c'est-\`{a}-dire $M\in\SLZ$,
on dit que $q$ et~$q'$ sont
(proprement) \'{e}quivalentes, et on note $q\sim q'$.  Deux formes quadratiques \'{e}quivalentes ont
donc m\^{e}me discriminant.  D'autre part, si l'une est primitive, l'autre aussi.  On montre que
l'ensemble $H(D)$ des classes de formes quadratiques primitives de discriminant~$D$ est fini.

  Les matrices $S=\smallmatrix0{-1}10$ et $T=\smallmatrix1101$ engendrent le groupe $\SLZ$.
Leur action (et plus g\'{e}n\'{e}ralement, l'action de $T^k=\smallmatrix1k01$ avec $k\in\Z$) sur les
formes est donn\'{e}e par
$$
  \eqalign{
     (a, b, c)\cdot S &= (c, -b, a),  \cr
     (a, b, c)\cdot T^k &= (a, b + 2ka, c + kb + k^2a). \cr}
$$
En appliquant certaines de ces transformations selon un algorithme simple, on obtient apr\`{e}s
un nombre fini d'\'{e}tapes une forme dite {\it r\'{e}duite} \'{e}quivalente \`{a} la forme initiale --- et
cet algorithme de r\'{e}duction est essentiel pour la th\'{e}orie et les applications des formes
quadratiques --- mais cela ne sera pas d\'{e}velopp\'{e} ici.

\label{nesteddiscs}
\section Discriminants embo\^{\i}t\'{e}s

  Les r\'{e}sultats mis en relief ici sont pour l'essentiel \'{e}parpill\'{e}s 
dans~\bib{Buell}\bib{BuchmannVollmer}\bib{Cox}, mais la notion de formes semi-\'{e}quivalentes 
que j'introduis ci-dessous n'y est pas pr\'{e}sente explicitement.
Nous consid\'{e}rons deux discriminants, $D$ et $Df^2$ o\`{u} $f$ est un entier~$>1$.

\subsection Matrices de remont\'{e}e

Pour $g$ entier tel que $0\leq g<f$, on pose $R_g=\smallmatrix fg01$.  On posera
aussi $R_f=\smallmatrix 100f$.

\label{nonequivMg}
\stat Lemme.  Les matrices $R_g$ sont de d\'{e}terminant $f$ et deux-\`{a}-deux non
\'{e}quivalentes sous l'action de~$\SLZ$.

\dem L'assertion sur le d\'{e}terminant est claire~; montrons la deuxi\`{e}me.  Supposons que $R_g$ et~$R_h$
sont \'{e}quivalentes, c'est-\`{a}-dire qu'il existe une matrice $\smallmatrix prst\in\SLZ$ telle que
$\smallmatrix fg01\smallmatrix prst=\smallmatrix fh01$.
\spot  Traitons d'abord le cas $g,h\neq f$.  On a alors
$$
  \pmatrix{f&g\cr0&1\cr}
  \pmatrix{p&r\cr s&t\cr}
  =
  \pmatrix{fp+gs & fr+gt\cr  s & t\cr}
  =
  \pmatrix{f&h\cr0&1\cr}.
$$
On obtient alors $s=0$, $t=1$, puis $p=1$ et $fr+g=h$.  Comme $0\leq g,h<f$ on a donc $g=h$.
\spot Dans le cas $g\neq f$ et $h=f$, on obtient $s=0$, $t=f$, puis $fp=1$ et $r+g=0$ ce qui est
impossible pour $f>1$.  \qed

\label{allequivMg}
\stat Lemme.  Soient $f>1$ un nombre premier, et $\smallmatrix prst$ une matrice \`{a} coefficients
entiers et de d\'{e}terminant~$f$.  Alors, il existe un unique $g\in\{0,\ldots,f\}$ et une matrice
$U\in\SLZ$ tels que $\smallmatrix prst=R_gU$.

\dem L'unicit\'{e} de~$g$ r\'{e}sultera du lemme~\ref{nonequivMg}.  Montrons l'existence.  Puisque
$pt-rs=f$, on a $\pgcd(s,t)\mid f$.  Comme $f$ est premier, ce pgcd vaut 1 ou~$f$.  

  Supposons d'abord que $\pgcd(s,t)=1$, et notons $\lambda$ et~$\mu$ des coefficients de B\'{e}zout~: 
$\lambda s+\mu t=1$.  On v\'{e}rifie alors que 
$$
  \pmatrix{p&r\cr s&t\cr} \pmatrix{t&\lambda\cr -s&\mu\cr}
  =
  \pmatrix{f&\lambda p+\mu r\cr 0&1\cr}.
$$
Soit $k$ entier tel que $0\leq\lambda p+\mu r-kf<f$.  On a ensuite
$$
  \pmatrix{f&\lambda p+\mu r\cr 0&1\cr} \pmatrix{1&-k\cr 0&1\cr}
  =
  \pmatrix{f&g\cr 0&1\cr}
  = R_g
  \qquad\hbox{o\`{u} $g=\lambda p+\mu r-kf$.}
$$ 
En posant $U^{-1}=\smallmatrix t\lambda{-s}\mu\smallmatrix 1{-k}01$, on a bien 
$\smallmatrix prstU^{-1}=R_g$.  De plus, $U=\smallmatrix{\mu+ks}{kt-\lambda}st\in\SLZ$.

  Il reste le cas o\`{u} $\pgcd(s,t)=f$.  On v\'{e}rifie que
$$
  \pmatrix{p&r\cr s&t\cr} \pmatrix{t/f&-r\cr -s/f&p\cr}
  =
  \pmatrix{1&0\cr 0&f\cr}.
$$
On a donc la relation cherch\'{e}e, avec $g=f$ et 
$U=\smallmatrix{t/f}{-r}{-s/f}p^{-1}=\smallmatrix pr{s/f}{t/f}$.  \qed

\subsection La remont\'{e}e principale

  Nous allons expliciter les liens entre les matrices de discriminant~$D$ et celles de
discriminant~$Df^2$.  Pour commencer, cette sous-section \'{e}nonce une proposition selon laquelle
toute forme de discriminant $Df^2$ est obtenue \`{a} partir d'une forme de discriminant~$D$, en
lui appliquant une matrice de remont\'{e}e sp\'{e}cifique (\`{a} savoir~$R_f$).

\label{Rfaction}
\stat Propri\'{e}t\'{e}.  Soient $f\in\N^*$ et $q=ax^2+bxy+cy^2$ une forme quadratique primitive de
discriminant~$D$ et telle que $\pgcd(a,f)=1$.  Alors $Q=q\cdot R_f$ (on a donc $Q(x,y)=q(x,fy)$)
est une forme quadratique primitive de discriminant~$Df^2$.

\dem La valeur du discriminant r\'{e}sulte de la proposition~\ref{discaction}.  Posons
$Q=Ax^2+Bxy+Cy^2$.  On a donc $A=a$, $B=fb$ et $C=f^2c$.  Tout facteur premier $p$ commun
\`{a} $A$, $B$ et~$C$ doit aussi diviser $b$ et $c$ puisque $\pgcd(a,f)=1$.  \qed

\label{normalform}
\stat Lemme.  Soit $Q=Ax^2+Bxy+Cy^2$ une forme quadratique de discriminant $Df^2$.  On suppose
que $\pgcd(2A,f)=1$.   Alors $Q$ est \'{e}quivalente \`{a} une forme quadratique $Q'=A'x^2+B'xy+C'y^2$
telle que $f\mid B'$ et $f^2\mid C'$ (et aussi $A'=A$).

\dem D'apr\`{e}s B\'{e}zout, il existe deux entiers $\lambda$ et~$\mu$ tels que $2\lambda A+\mu f=1$.
Posons $M=\smallmatrix1{-\lambda B}01$ et $Q'=Q\cdot M$.  On a alors $A'=A$ et 
$$
  B' =
  B-2A\lambda B
  =
  B - (1 - \mu f)B
  \equiv 0 
  \pmod f.
$$
De plus, on a $\det M=1$ donc $\disc Q'=\disc Q=Df^2$.  Donc $f^2\mid B'^2-4AC'$.  Puisque 
$f\mid B'$ et $\pgcd(f,4A)=1$, on obtient $f^2\mid C'$.  \qed

\label{Qfromq}
\stat Proposition.  On suppose $f$ impair.  Soit $Q$ une forme primitive de discriminant~$Df^2$.
Il existe une forme primitive $q$ de discriminant~$D$ telle que $Q=q\cdot R_fU$ o\`{u} $U\in\SLZ$.

\dem Posons $Q=Ax^2+Bxy+Cy^2$.  Le lemme~\ref{normprimetof}, permet de supposer que~$\pgcd(A,f)=1$.
D'apr\`{e}s le lemme~\ref{normalform}, on peut supposer que $f\mid B$ et $f^2\mid C$.  Il suffit alors
de poser $q=Q\cdot\smallmatrix100{1/f}$.  \qed

\subsection Formes semi-\'{e}quivalentes

  J'introduis ici la notion de formes semi-\'{e}quivalentes.  Le th\'{e}or\`{e}me qui suivra permet de
renforcer les propri\'{e}t\'{e}s des formes semi-\'{e}quivalentes, lorsque $f$ est premier.

\label{defsemiequiv}
\statrm D\'{e}finition.  Soient $Q_1$ et $Q_2$ deux formes primitives de discriminant~$Df^2$.  Je
dis que $Q_1$ et~$Q_2$ sont {\it semi-\'{e}quivalentes} (ou {\it fondamentalement \'{e}quivalentes} lorsque
le discriminant~$D$ est fondamental) si il existe deux formes primitives \'{e}quivalentes
$q_1$ et $q_2$ de discriminant $D$, deux entiers $g_1$ et $g_2$ avec $0\leq g_1,g_2\leq f$, ainsi
que deux $\SLZ$-matrices $U_1,U_2$, tels que $Q_i=q_i\cdot R_{g_i}U_i$ (pour $i=1,2$).

  La semi-\'{e}quivalence est une relation d'\'{e}quivalence moins fine que l'\'{e}quivalence : deux formes
\'{e}quivalentes sont n\'{e}cessairement semi-\'{e}quivalentes.

\label{Aurore}%label{RgRf}
\stat Lemme.  Soient $D$ un discriminant et $f\in\N^*$.
Soient $q$ et~$Q$ deux formes primitives de discriminants respectifs $D$ et~$Df^2$
telles que $Q=q\cdot R_gU$ pour en entier $h$ ($0\leq h\leq f$) et une $\SLZ$-matrice $U$.
Il existe une forme $q_0$ \'{e}quivalente \`{a}~$q$ et $V\in\SLZ$ telles que $Q=q_0\cdot R_fV$.

\dem Posons $q=ax^2+bxy+cy^2$.  Si $g=f$, c'est clair.  Nous supposerons donc que $g<f$.  Dans ce
cas, $q\cdot R_g$ est alors donn\'{e}e par
$$
  Q_g = q \cdot R_g = Ax^2 + BXY + CY^2
  \qquad\hbox{avec $A=f^2a$, $B=f(2ag+b)$ et $C=ag^2+bg+c$.}
  \eqdef{qRh}
$$
Appliquons $S=\smallmatrix0{-1}10$ puis $R_f^{-1}=\smallmatrix100{1/f}$~:
$$
  \eqalign{   
  Q_g\cdot SR_f^{-1}
  &= Cx^2 - {B\over f}xy + {A\over f^2}y^2
   = (ag^2 + bg + c) x^2 - (2ag + b) xy + a y^2 \cr
  &= a(y - gx)^2 - b(y - gx)x + cx^2 \cr
  &= (ax^2+bxy+cy^2) \cdot  \pmatrix{g&-1\cr1&0\cr}
   = q \cdot\pmatrix{g&-1\cr1&0\cr}. \cr}
$$
Posons $q_0=q\cdot\smallmatrix g{-1}10$.  On a alors $Q_g=q_0\cdot R_fS^{-1}$, donc 
$Q=q_0\cdot R_f(S^{-1}U)$.
\qed
  
  Lorsque deux formes sont semi-\'{e}quivalentes, et que $f$ est premier, on peut affirmer que deux
autres conditions, similaires mais un peu plus contraignantes que celle de la d\'{e}finition,
sont satisfaites~:

\label{semiequiv}
\stat Th\'{e}or\`{e}me.  On suppose $f$ premier.  Soient $Q_1$ et $Q_2$ deux formes primitives de
discriminant~$Df^2$.  Alors chacune des deux conditions suivantes est satisfaite si et seulement
si $Q_1$ et $Q_2$ sont semi-\'{e}quivalentes.
\item{(a)} Il existe une forme primitive~$q$ de discriminant~$D$, deux entiers $g_1$ et~$g_2$
avec $0\leq g_1,g_2\leq f$, et deux $\SLZ$-matrices $U_1,U_2$ tels que $Q_i=q\cdot R_{g_i}U_i$
(pour $i=1,2$).
\item{(b)} Il existe deux formes primitives \'{e}quivalentes $q_1$ et~$q_2$ de discriminant~$D$ et deux 
$\SLZ$-matrices $V_1,V_2$ telles que $Q_i=q_i\cdot R_fV_i$ (pour $i=1,2$).

\dem Les conditions (a) et (b) sont {\sl a priori} plus contraignantes que celle utilis\'ee
dans~\ref{defsemiequiv} pour d\'efinir la notion de formes semi-\'{e}quivalentes.  Il s'agit donc de
montrer que, si $Q_1$ et $Q_2$ sont semi-\'{e}quivalentes,
elles satisfont (a) et (b).  Supposons donc $Q_1$ et $Q_2$ semi-\'{e}quivalentes, et reprenons les
notations $q_1$, $q_2$, $g_1$, $g_2$, $U_1$, $U_2$ de la d\'{e}finition~\ref{defsemiequiv}.

  Soit $U\in\SLZ$ telle que $q_1\cdot U=q_2$.  On a $Q_2=q_1\cdot UR_{g_2}U_2$.
D'apr\`{e}s le lemme~\ref{allequivMg}, on peut \'{e}crire $UR_{g_2}=R_gU'$ avec $0\leq g\leq f$ et 
$U'\in\SLZ$.  On obtient donc (a) avec $q=q_1$, en utilisant $g$ dans le r\^{o}le de~$g_2$, et $U'U_2$
dans celui de $U_2$. 

  D'autre part, le lemme~\ref{Aurore} indique qu'il existe une forme $q'_1$ \'{e}quivalente \`{a}~$q_1$
et une matrice $V_1\in\SLZ$ telles que $Q_1=q'_1R_fV_1$.   De m\^{e}me, on a $Q_2=q'_2R_fV_2$ pour une
forme $q'_2$ \'{e}quivalente \`{a}~$q_2$ et $V_2\in\SLZ$.  On a obtenu (b) puisque
$q'_1\sim q_1\sim q_2\sim q_2'$.  \qed

\subsection La descente

  On verra en particulier ici que si deux formes semi-\'{e}quivalentes sont obtenues \`{a} partir
de deux formes $q_1$ et~$q_2$ de discriminant~$D$, alors $q_1\sim q_2$.

\label{integralRf}
\stat Lemme.  Soient $U\in\SLZ$ et $V=R_fUR_f^{-1}$.  Soit $q=ax^2+bxy+cy^2$ une forme telle que
$\pgcd(a,f)=1$.  On suppose que $q\cdot V$ est une forme, c'est-\`{a}-dire que ses coefficients sont
entiers.  Alors $V\in\SLZ$.

\dem On a $\det V=1$.  Il s'agit donc de montrer que la matrice~$V$ est \`{a} coefficients entiers.
Posons $U=\smallmatrix prst$.  On a donc $V=\smallmatrix p{r/f}{fs}t$.   Posons
$q\cdot V=a'x^2+b'xy+c'y^2$.  En vertu de~\eqref{coeffsequivform}, on a 
$c' = q(r/f, t)=a(r/f)^2+brt/f+ct^2$.  L'hypoth\`{e}se indique que $c'$ est entier, donc 
$f^2\mid ar^2+brtf=r(ar+btf)$.  Soit $p^k$ une puissance maximale de $p$ premier telle que 
$p^k\mid f$ (avec $k\geq1$).    Montrons par r\'{e}currence que $p^k\mid r$.  Si on a montr\'{e} que
$p^l\mid r$ (avec $l<k$), alors on peut \'{e}crire 
$$
  p \Bigmid p^{2(l-1)} \Bigmid {r\over p^l} \bigg(a {r\over p^l} + bt {f\over p^l}\bigg).
$$
Sachant que $p\notmid a$, on obtient ais\'{e}ment que $p\mid r/p^l$.   La r\'{e}currence montre donc que 
$p^k\mid r$.  Finalement $f\mid r$.  \qed

\label{semiequivequiv}
\stat Proposition.  On suppose $f$ premier.  Soient $Q_1$ et~$Q_2$ des formes primitives de
discriminant $Df^2$ semi-\'{e}quivalentes.  Soient $q_1,q_2$ deux formes primitives de
discriminant~$D$ et~$M_1,M_2$ deux matrices enti\`{e}res de d\'{e}terminant~$f$ telles 
que~$Q_i=q_i\cdot M_i$ (pour $i=1,2$). (Les formes $q_i$ et les matrices $M_i$ existent d'apr\`{e}s la
proposition~\ref{Qfromq}.)  Alors $q_1\sim q_2$.

\dem Tout d'abord, en utilisant le lemme~\ref{normprimetof}, on peut supposer les coefficients
$a_1$ et~$a_2$ du terme en~$x^2$ de $q_1$ et~$q_2$ v\'{e}rifient $\pgcd(a_i,f)=1$.
D'apr\`{e}s le (b) du th\'{e}or\`{e}me~\ref{semiequiv}, il existe deux formes \'{e}quivalentes $q'_1$ et
$q'_2$ telles que $Q_i=q'_i\cdot R_fV_i$ avec $V_i\in\SLZ$.  On a alors 
$q_i\cdot M_iV_i^{-1}R_f^{-1}=q'_i$.  Le lemme~\ref{allequivMg} indique que $M_i$ s'\'{e}crit sous la
forme $R_{g_i}U_i$, puis le lemme~\ref{Aurore} permet d'\'{e}crire 
$q''_i\cdot R_fV'_iV_i^{-1}R_f^{-1}=q'_i$ o\`{u} $q''_i\sim q_i$ et $V'_i\in\SLZ$.  Le
lemme~\ref{integralRf} pr\'{e}cise alors que le produit $R_fV'_iV_i^{-1}R_f^{-1}$ est \'{e}l\'{e}ment de~$\SLZ$.
Cela signifie que $q''_i\sim q_i'$ donc $q_i\sim q'_i$.  Comme $q'_1\sim q'_2$, on obtient 
$q_1\sim q_2$.  \qed

\label{drown}
\stat Th\'{e}or\`{e}me.  Soient $D$ un discriminant et $f$ un nombre premier impair.  Pour $Q$ forme
quadratique primitive de discriminant~$Df^2$, la proposition~\ref{Qfromq} pr\'{e}cise qu'il existe une
forme quadratique $q$ de discriminant~$D$ telle que $q\cdot R_f\sim Q$.  La correspondance
$Q\mapsto q$ ainsi obtenue d\'{e}finit une application surjective $\pi$ de l'ensemble des classes
de formes quadratiques primitives de discriminant $Df^2$ sur celui des formes quadratiques
primitives de discriminant~$D$.

\dem La proposition~\ref{semiequivequiv} montre que la forme $q$ obtenue ne d\'{e}pend pas de la
repr\'{e}sentante~$Q$ choisie au sein d'une m\^{e}me classe d'\'{e}quivalence (m\^{e}me de semi-\'{e}quivalence).  Donc
on obtient bien une application.  De plus, pour $q$ primitive de discriminant~$D$, la
proposition~\ref{Rfaction} fournit un ant\'{e}c\'{e}dent (la classe de $q\cdot R_f$) \`{a} la classe de~$q$. 
\qed

\subsection Les autres remont\'{e}es

  Enfin, ici on d\'{e}nombre et explicite les formes et classes de formes de discriminant $Df^2$
obtenues \`{a} partir d'une forme de discriminant~$D$.

% Buell, page 113
% Aurore, page 92
\label{qRgprim}
\stat Proposition.  On suppose que $f$ est un nombre premier, que $q=ax^2+bxy+cy^2$ est
primitive et que $\pgcd(a,f)=1$.  Alors $Q_g=q\cdot R_g$ est primitive pour exactement 
$f-\Legendre Df$ valeurs de $g$ (celles telles que $q(g,1)\not\equiv0$ modulo~$f$ ; ici
$\Legendre Df$ est un symbole de Kronecker).

\dem On a d\'{e}j\`{a} vu (proposition~\ref{Rfaction}) que $Q_f$ est primitive.  Supposons donc $g\neq f$. 
Dans ce cas, $Q_g$ et ses coefficients $A,B,C$ sont donn\'{e}s par la formule~\eqref{qRh}.  Puisque
$a,b,c$ sont globalement premiers entre eux, les entiers $a$, $2ag+b$ et $ag^2+bg+c$ le sont aussi.
Donc $Q_g$ est primitive si et seulement si $\pgcd(f,C)=1$.  Mais $C=q(g,1)$ est un
polyn\^{o}me en $g$ de discriminant~$D$.  Il s'annule modulo~$f$ pour $1+\Legendre Df$ valeurs de~$g$.
Au total $1+f-(1+\Legendre Df)=f-\Legendre Df$ valeurs de~$g$ donnent une forme $Q_g$ primitive.
\qed

\label{nonequivQg}
\stat Proposition.  Soient $q$ une forme de discriminant~$D$ n\'{e}gatif et $f$ un nombre premier
impair. Notons $Q_g=q\cdot R_g$ pour $0\leq g\leq f$ et $q(g,1)\not\equiv0$ modulo~$f$.  Les
$f-\Legendre Df$ formes primitives~$Q_g$ obtenues sont deux \`{a} deux non-\'{e}quivalentes, sauf pour les
discriminants exceptionnels $D=-3$ et $D=-4$ pour lesquels respectivement exactement 2 et 3 valeurs
distinctes de~$g$ donnent la m\^{e}me classe d'\'{e}quivalence.

\dem Supposons $Q_{g_1}$ et~$Q_{g_2}$ \'{e}quivalentes, et d\'{e}signons par~$N$ une $\SLZ$-matrice telle
que $Q_{g_1}\cdot N=Q_{g_2}$.    On a alors
$$
  q \cdot R_{g_1} N R_{g_2}^{-1} = q
  \eqdef{qtoq}
$$

  Montrons d'abord que $R_{g_1}NR_{g_2}^{-1}\in\SLZ$.  D'apr\`{e}s le lemme~\ref{Aurore}, il existe deux
formes $q_1$ et~$q_2$ \'{e}quivalentes \`{a}~$q$ ainsi que deux $\SLZ$-matrices $V_1$ et~$V_2$ telles que 
$Q_{g_i}=q_i\cdot R_fV_i$ (pour $i=1,2$).   En d\'{e}signant $G_i$ deux $\SLZ$-matrices telles que
$q_i=q\cdot G_i$ (pour $i=1,2$), on a
$$
  R_{g_1} N R_{g_2}^{-1} = G_1 R_f V_1 N V_2^{-1} R_f^{-1} G_2^{-1}.
$$
Mais le lemme~\ref{integralRf} pr\'{e}cise que $R_f(V_1NV_2^{-1})R_f^{-1}$ est \'{e}l\'{e}ment de~$\SLZ$, donc 
$R_{g_1}NR_{g_2}^{-1}$ aussi.

  L'\'{e}quation~\eqref{qtoq} indique donc que $R_{g_1}NR_{g_2}^{-1}$ est un automorphisme~$Z$ de~$q$.
Si $D<-4$, alors $Z=\pm I$.  On en d\'{e}duit $R_{g_1}(\pm N) = R_{g_2}$, ce qui implique $g_1=g_2$
d'apr\`{e}s le lemme~\ref{nonequivMg}.

  Si $D=-4$, il y a un cas alternatif qui est $Z=\pm\smallmatrix{-b/2}{-c}a{b/2}$.  On peut
simplifier en prenant $q(x,y)=x^2+y^2$ et donc $Z=\pm S$.  On obtient $R_{g_1}N=\pm SR_{g_2}$.  Le
lemme~\ref{allequivMg} permet d'\'{e}crire $SR_{g_2}=\smallmatrix0{-1}f{g_2}=R_{g'}M$ avec
$g'=-g_2^{-1}\bmod f$ lorsque $0<g_2<f$ (ou $g'=f$ lorsque $g_2=0$~; ou $g'=0$ lorsque $g_2=f$)
et $M\in\SLZ$.  Enfin, le lemme~\ref{nonequivMg} pr\'{e}cise que $g_1=g'$.  On remarquera que
$q(g_i,1)=g_i^2+1\not\equiv0$ modulo~$f$ implique $g_1\neq g_2$.  On obtient donc exactement deux
fois chaque classe.

  Si $D=-3$, il y a deux cas alternatifs.  On simplifie en prenant $q(x,y)=x^2+xy+y^2$ et
$Z=\pm\smallmatrix0{-1}11$ ou $Z=\pm\smallmatrix11{-1}0$.  On trouve 
$Q_g\sim Q_{g'}\sim Q_{g''}$ lorsque $g'\equiv-g^{-1}-1$ et $g''\equiv-(g+1)^{-1}$ ou lorsque
$\{g,g',g''\}=\{0,f-1,f\}$.  L\`{a} aussi, $q(g,1)=g^2+g+1\neq0$ implique que $g$, $g'$ et~$g''$ sont
distincts.  \qed

\label{allequivQg}
\stat Proposition.  On suppose que $f$ est un nombre premier impair.  Soit $Q$ une forme
quadratique primitive de discriminant $Df^2$ et $q\in\pi(Q)$.  Alors $Q$ est \'{e}quivalente \`{a} l'une
des $q\cdot R_g$ pour un $g$ tel que $f\notmid q(g,1)$. 

\dem Soit $R$ une matrice de d\'{e}terminant $f$ telle que $q\cdot R=Q$ (cette
matrice existe puisque $q=\pi(Q)$).  Soient $g$ et~$M$ (donn\'{e}s par le lemme~\ref{allequivMg}),
tels que $R=R_gM$.  On a donc $Q=q\cdot R_gM\sim q\cdot R_g$.  De plus, $Q\cdot R_g$ est primitive
donc $f\notmid q(g,1)$, d'apr\`{e}s~\ref{qRgprim}.  \qed

\section Bibliographie 

\showbib

\bye